\documentclass[a4paper,12pt]{article}
\usepackage{amsmath}
\usepackage{amssymb}
\def\cT{\mathcal{T}}
\def\cN{\mathcal{N}}
\def\cM{\mathcal{M}}
\newcommand{\be}{\begin{equation}}
\newcommand{\ee}{\end{equation}}
\newcommand{\ra}{\rightarrow}

\newcommand{\bea}{\begin{eqnarray}}
\newcommand{\eea}{\end{eqnarray}}
\newcommand{\beas}{\begin{eqnarray*}}
\newcommand{\eeas}{\end{eqnarray*}}
\newcommand{\R}{\mathbb{R}}

\newcommand{\nn}{\nonumber}
\newcommand{\ot}{\otimes}
\newcommand{\op}{\oplus}
\newcommand{\wt}{\widetilde}
\newcommand{\s}{\circ}
\def\Sec{\text{Sec}}
\def\sT{\text{T}}
\def\xd{\text{d}}
\newcommand{\pa}{\partial}

\newcommand{\A}{{\cal A}}
\def\la{\langle}
\def\ran{\rangle}
\mathchardef\za="710B  
\mathchardef\zb="710C  
\mathchardef\zg="710D  
\mathchardef\zd="710E  
\mathchardef\zve="710F 
\mathchardef\zz="7110  
\mathchardef\zh="7111  
\mathchardef\zvy="7112 
\mathchardef\zi="7113  
\mathchardef\zk="7114  
\mathchardef\zl="7115  
\mathchardef\zm="7116  
\mathchardef\zn="7117  
\mathchardef\zx="7118  
\mathchardef\zp="7119  
\mathchardef\zr="711A  
\mathchardef\zs="711B  
\mathchardef\zt="711C  
\mathchardef\zu="711D  
\mathchardef\zvf="711E 
\mathchardef\zq="711F  
\mathchardef\zc="7120  
\mathchardef\zw="7121  
\mathchardef\ze="7122  
\mathchardef\zy="7123  
\mathchardef\zf="7124  
\mathchardef\zvr="7125 
\mathchardef\zvs="7126 
\mathchardef\zf="7127  
\mathchardef\zG="7000  
\mathchardef\zD="7001  
\mathchardef\zY="7002  
\mathchardef\zL="7003  
\mathchardef\zX="7004  
\mathchardef\zP="7005  
\mathchardef\zS="7006  
\mathchardef\zU="7007  
\mathchardef\zF="7008  
\mathchardef\zW="700A  

\newcommand{\epf}{\hfill$\Box$}
\newcommand{\bepf}{\textit{Proof.-} }

\begin{document}

\title{Courant-Nijenhuis tensors and
generalized geometries}
\author{
Janusz Grabowski\thanks{Supported by KBN, grant No 2 P03A 020 24.}\\
Mathematical Institute, Polish Academy of Sciences\\ ul.
\'Sniadeckich 8, P. O. Box 21, 00-956 Warszawa, Poland\\ {\it
e-mail:} jagrab@impan.gov.pl} \date{}\maketitle

\newtheorem{re}{Remark}
\newtheorem{theo}{Theorem}
\newtheorem{prop}{Proposition}
\newtheorem{lem}{Lemma}
\newtheorem{cor}{Corollary}
\newtheorem{ex}{Example}

\rightline{\it Dedicated to Jos\'e~F.~Cari\~nena on his 60th
birthday}

\bigskip
\begin{abstract} Nijenhuis tensors $N$ on Courant algebroids
compatible with the pairing are studied. This compatibility
condition turns out to be of the form $N+N^*=\zl I$ for
irreducible Courant algebroids, in particular for the extended
tangent bundles $\cT M=\sT M\op \sT^* M$. It is proved that
compatible Nijenhuis tensors on irreducible Courant algebroids
must satisfy quadratic relations $N^2-\zl N+\zg I=0$, so that the
corresponding hierarchy is very poor. The particular case $N^2=-I$
is associated with Hitchin's generalized geometries and the cases
$N^2=I$ and $N^2=0$ -- to other "generalized geometries". These
concepts find a natural description in terms of supersymplectic
Poisson brackets on graded supermanifolds.

\bigskip\noindent
\textit{\textbf{MSC 2000:} Primary 17B99; Secondary 17B62, 53C15,
53C56, 53D05, 53D17.}

\medskip\noindent
\textit{\textbf{Key words:} Nijenhuis tensor, Leibniz algebra,
Courant bracket, Courant algebroid, generalized complex geometry.}

\end{abstract}

\section{Introduction}
The theory of Nijenhuis tensors on Lie algebras goes back to a
concept of contractions of Lie algebras introduced by
E.~J.~Saletan \cite{Sa}. The study of Nijenhuis tensors for Lie
algebroids and Nijenhuis tensors on Poisson manifolds have been
originated in \cite{MM,KSM}. In \cite{CGM0} it has been developed
the theory of Nijenhuis tensors for associative products, and in
\cite{CGM} contractions and Nijenhuis tensors have been studied
for algebraic operations of arbitrary type on sections of vector
bundles.

Recall that a \textit{Nijenhuis tensor} $N$ for a bilinear
operation "$\s$" on sections of a vector bundle $A$ over $M$ is a
$(1,1)$-tensor $N\in\Sec(A\ot A^*)$, viewed also as vector bundle
morphism $N:A\ra A$ (or the corresponding $C^\infty(M)$-linear map
$N:\Sec(A)\ra\Sec(A)$ on sections), such that its
\textit{Nijenhuis torsion}
\begin{equation}\label{1}
\text{Tor}_N(X,Y)=N(X)\s N(Y)-N(X\s_NY)
\end{equation}
vanishes. Here "$\s_N$" is the `contracted' product:
\begin{equation}\label{2}
X\s_NY=N(X)\s Y+X\s N(Y)-N(X\s Y).
\end{equation}

This general procedure has been applied in \cite{CGMc} to Leibniz
algebras and \textit{Courant algebroids} \cite{LWX} in their
Leibniz algebra formulation. {\it Leibniz algebras} --
non-skew-symmetric generalizations of Lie algebras -- were studied
first by J.-L.~Loday \cite{Lo} (they are called sometimes
\textit{Loday algebras}) and the (co)homology theory of Lie
algebras was generalized to this framework.

\medskip\noindent
\textbf{Definition 1} A \textit{Leibniz product (bracket)} on a
vector space $\A$ is a bilinear operation "$\s$" satisfying the
Jacobi identity
\begin{equation}\label{3}
(X\s Y)\s Z=X\s(Y\s Z)-Y\s(X\s Z)
\end{equation}
for all $X,Y,Z\in\A$. The space $\A$ with a Leibniz product we
call a \textit{Leibniz algebra}.

\medskip\noindent
Let now "$\s$" be a local Leibniz product on the space $\Sec(A)$
of sections of a vector bundle $A$ over $M$, i.e. a product which
is locally defined by a bidifferential operator, and let $N:A\ra
A$ be a $(1,1)$-tensor over $A$. According to the general scheme
in \cite{CGM}, if the Nijenhuis torsion (\ref{1}) vanishes, the
contracted product (\ref{2}) is a Leibniz product which is
\textit{compatible} with the original one, i.e. $X\s_N Y+\zl X\s
Y$ is a Leibniz product for any $\zl\in\R$. Note that the
compatibility is always satisfied.

\begin{theo}\label{t0}\cite{CGMc}
The products "$\s_N$" and "$\s$" are always compatible.
The contracted product (\ref{2}) is still Leibniz if and only if
the Nijenhuis torsion (\ref{1}) is a 2-cocycle with respect to the
Leibniz cohomology operator, i.e.
\bea\label{coc}
(\zd \text{Tor}_N)(X,Y,Z)&=&\text{Tor}_N(X,Y\s Z)-\text{Tor}_N(X\s Y,Z)-\text{Tor}_N(Y,X\s Z)\\
&&-\text{Tor}_N(X,Y)\s Z+X\s \text{Tor}_N(Y,Z)-Y\s
\text{Tor}_N(X,Z)=0.\nn
\eea
In this case "$\s_N$" and "$\s$" are compatible Leibniz products.
\end{theo}

\medskip\noindent
The tensor $N$ we will call a \textit{Nijenhuis tensor} (for the
Leibniz algebra $\Sec(A)$) if the Nijenhuis torsion $\text{Tor}_N$
vanishes and a \textit{weak Nijenhuis tensor} if the Nijenhuis
torsion $\text{Tor}_N$ is a Leibniz 2-cocycle. In both cases the
contracted product "$\s_N$" is Leibniz and it is compatible with
the original one.

\medskip\noindent
\textbf{Example 1} An interesting example of a Leibniz product is
the following Leibniz algebra version of the \textit{Courant
bracket} (called sometimes also \textit{Dorfman bracket}) on
sections $X+\zx$ of the bundle $\cT M=\sT M\op \sT^*M$:
\begin{equation}\label{Cou}(X+\zx)\s(Y+\zh)=[X,Y]+(\pounds_X\zh-i_Yd\zx).
\end{equation}
Here $[X,Y]$ is clearly the bracket of vector fields, $\pounds_X$
is the Lie derivative, etc. The extended tangent bundle $\cT M$
with the canonical symmetric pairing, coming from the contraction,
and with the Courant bracket is an example of a \textit{Courant
algebroid} (cf. \cite{LWX,Ro}).

\medskip\noindent
Since a Courant algebroid (see the next section) is not only a
Leibniz algebra on sections of a vector bundle but also a
non-degenerate pairing with certain consistency conditions with
the Leibniz product, it has been studied in \cite{CGMc} what is
the property of $N$ that ensures the consistency conditions being
satisfied also for "$\s_N$". It turns out that it is sufficient to
assume that $N+N^*=\zl I$, $\zl\in\R$, where $N^*$ is dual to $N$
with respect to the pairing. This implies in the particular case
of $\cT M$ that such $N$ is associated with a triplet consisting
of a $(1,1)$-tensor, a 2-form, and a bivector field on $M$.

In this paper we prove that this condition is also necessary for
so called \textit{irreducible} Courant algebroids ($\cT M$ is a
canonical example). We prove also that such compatible Nijenhuis
tensors on irreducible Courant algebroids must satisfy
additionally a quadratic equation $N^2-\zl N+\zg I=0$, se the
associated hierarchy is trivial. Particular cases: $N^2=-I$,
$N^2=I$, and $N^2=0$ correspond to the so called \textit{complex},
\textit{product}, and \textit{tangent Courant structures},
respectively. The complex Courant structures on $\cT M$ were
introduced recently by N.~Hitchin \cite{Hi} under the name of
\textit{complex generalized geometries} and they drew much
attention among mathematicians and physicists. Our work shows
that, in practice, due to the above quadratic equation, no more
"generalized geometries" in this sense than complex, product, and
tangent are possible. Since, according to \cite{Ro2}, any Courant
algebroid is associated with a cubic homological Hamiltonian $\zY$
on a symplectic $\cN$-manifold of degree 2, we show that in this
language complex Courant structures correspond to certain
quadratic super-functions $N$ such that $\{\{\zY,N\},N\}=-\zY$,
where the bracket is the corresponding Poisson superbracket.

\section{Nijenhuis tensors for Courant algebroids}
Let us recall briefly the structure of a Courant algebroid. We
will use here the Leibniz product (bracket) version of the Courant
bracket presented already in \cite{Ro} with some simplifications
discussed already in \cite{CGMc} (cf. also \cite[Definition1]{GM},
\cite[Definition 2.1]{KS1}, and \cite{Uch}).

\medskip\noindent
\textbf{Definition 2\ } A \textit{Courant algebroid} is a vector
bundle $\zt:A\ra M$ equipped with a Leibniz product (bracket)
"$\s$" on $Sec(A)$, a vector bundle map (covering the identity)
$\zr:A\ra \sT M$ and a nondegenerate symmetric bilinear form
$\la{\cdot},{\cdot}\ran$ on $A$ satisfying the identities
\bea\label{4}
&\zr(X)\la Y,Y\ran=2\la X,Y\s Y\ran,\\ &\zr(X)\la Y,Y\ran=2\la X\s
Y,Y\ran.\label{5}
\eea
Note that (\ref{4}) is equivalent to
\begin{equation}\label{4a}
\zr(X)\la Y,Z\ran=\la X,Y\s Z+Z\s Y\ran.
\end{equation}
Similarly, (\ref{5}) easily implies the invariance of the pairing
$\la{\cdot},{\cdot}\ran$ with respect to the left multiplication
\begin{equation}\label{6}\zr(X)\la
Y,Z\ran=\la X\s Y,Z\ran+\la Y,X\s Z\ran
\end{equation}
and that $\zr$ is the anchor map for the left multiplication:
\begin{equation}\label{7}
X\s(fY)=fX\s Y+\zr(X)(f)Y.
\end{equation}
A rather unpleasant constatation is that, even when the Nijenhuis
torsion $\text{Tor}_N$ of a $(1,1)$-tensor $N\in\Sec(A^*\ot{A})$
vanishes (so the contracted bracket is a Leibniz bracket), the
conditions (\ref{4}) and (\ref{5}) need not to be satisfied
automatically for the `contracted' product (\ref{2}). Assume
therefore that $N$ is just a $(1,1)$-tensor on $A$ (do not assume
that $N$ is Nijenhuis at the moment) and repeat in short from
\cite{CGMc} the checking under what conditions the identities
(\ref{4}) and (\ref{5}) are still satisfied for "$\s_N$". Exactly
as in the classical case of a Lie algebroid contraction
\cite[Lemma 2]{CGM}, we have the anchor $\zr_N=\zr\circ N$ for the
contracted multiplication
\begin{equation}\label{8}
X\s_N(fY)=f(X\s_NY)+\zr(NX)(f)Y.
\end{equation}
Let $N^*$ be the adjoint of $N$ with respect to the pairing: $\la
NX,Y\ran=\la X,N^*Y\ran$ and let $\zD=N+N^*$. Using the invariance
(\ref{5}) we get easily
\beas
\la X\s_NY,Z\ran&=&\la NX\s Y+X\s NY-N(X\s Y),Z\ran\\
&=&\zr(NX)\la Y,Z\ran-\la Y,NX\s Z\ran+ \la Y,N^*(X\s Z)\ran+\la
Y, X\s N^*Z\ran,
\eeas
which equals $\zr(NX)\la Y,Z\ran-\la Y,X\s_NZ\ran$ if and only if
$\la Y,X\s \zD Z-\zD(X\s Z)\ran=0$ for all $X,Y,Z$, i.e. if and
only if $\zD$ commutes with the left multiplication
\begin{equation}\label{9}
X\s \zD Z-\zD(X\s Z)=0.
\end{equation}
Thus (\ref{9}) is equivalent to the invariance of the pairing with
respect to "$\s_N$":
$$\zr_N(X)\la Y,Z\ran=\la X\s_NY,Z\ran+\la Y,X\s_NZ\ran.
$$
Similarly, checking (\ref{4}) for "$\s_N$", we get
$$
\la X,Y\s_NY\ran=\frac{1}{2}\zr(X)\la Y,\zD
Y\ran-\frac{1}{2}\zr(N^*X)\la Y,Y\ran
$$
which equals $\frac{1}{2}\zr(NX)\la Y,Y\ran$ if and only if
$\zr(X)\la Y,\zD Y\ran=\zr(\zD X)\la Y,Y\ran$. The latter can be
rewritten in the form
$$\la X,Y\s \zD Y+\zD Y\s Y\ran=2\la\zD X,Y\s Y\ran$$
or
$$Y\s \zD Y+\zD Y\s Y=2\zD(Y\s Y).$$
Using (\ref{9}) we get finally the condition
\begin{equation}\label{10}\zD(Y\s Y)=\zD Y\s Y.
\end{equation}

\begin{theo} [\cite{CGMc}]\label{t1.1} If $N:A\ra A$ is a $(1,1)$-tensor on a
Courant algebroid, then the contracted product (\ref{2}) is
compatible with the symmetric pairing $\la{\cdot},{\cdot}\ran$ of
the Courant algebroid, in the sense that (\ref{4}) and (\ref{5})
are satisfied for "$\s_N$" and $\zr_N$, if and only if
$$X\s (N+N^*)
Y=(N+N^*)(X\s Y)\quad \text{and} \quad (N+N^*)(Y\s Y)=(N+N^*)Y\s
Y$$ for all sections $X,Y$ of $A$.
\end{theo}
It is clear that, how restrictive the above conditions are,
depends on `irreducibility' of the Courant product. However, there
is one (and only one) case which works for any Courant algebroid,
namely the case $N+N^*=\zl I$, $\zl\in\R$. A Courant algebroid we
call \textit{irreducible} if $\zl I$ are the only $(1,1)$-tensors
$\zD:A\ra A$ satisfying (\ref{9}) and (\ref{10}).

\begin{theo}\label{t1.2} The classical Courant algebroid structure on
$\cT M=\sT M\ot \sT^*M$ is irreducible.
\end{theo}
\bepf Suppose that the $(1,1)$-tensor $\zD$ commutes with the left
multiplication. In local coordinates $(x^i)$ we can write
$\zD(\pa_j)=\sum_i\left(\zD_j^i(x)\pa_i+\zD_j^{*i}(x)\xd
x^i\right)$ and $\zD(\xd x^j)=\sum_i\left(\zD_{*j}^i(x)\pa_i+
\zD_{*j}^{*i}(x)\xd x^i\right)$. In view of
$$0=\zD(\pa_k\s\pa_j)=\pa_k\s\zD(\pa_j)=\sum_i\left(\frac{\pa\zD^i_j}
{\pa x^k}(x)\pa_i+\frac{\pa\zD^{*i}_j}{\pa x^k}(x)\xd x^i\right)$$
and
$$0=\zD(\pa_k\s\xd x^j)=\pa_k\s\zD(\xd x^j)=\sum_i\left(\frac{\pa\zD^i_{*j}}
{\pa x^k}(x)\pa_i+\frac{\pa\zD^{*i}_{*j}}{\pa x^k}(x)\xd
x^i\right)$$ we get that $\zD^i_j(x)=\zD^i_j$,
$\zD^i_{j*}(x)=\zD^i_{*j}$, $\zD^{*i}_j(x)=\zD^{*i}_j$, and
$\zD^{*i}_{*j}(x)=\zD^{*i}_{*j}$ are constant. Now, since
$$(x^k\pa_j)\s\pa_k=[x^k\pa_j,\pa_k]=-\pa_j\quad \text{and}\quad (x^j\pa_k)\s\xd
x^k=\pounds_{x^j\pa_k}\xd x^k=\xd x^j,$$ we have
\be\label{1.r1}-\zD(\pa_j)=-\sum_i\left(\zD_j^i\pa_i+\zD_j^{*i}\xd
x^i\right)=x^k\pa_j\s\sum_i\left(\zD_k^i\pa_i+\zD_k^{*i}\xd
x^i\right)=-\zD^k_k\pa_j+\zD_k^{*j}\xd x^k\ee and
\be\label{1.r2}\zD(\xd x^j)=\sum_i\left(\zD_{*j}^i\pa_i+\zD_{*j}^{*i}\xd
x^i\right)=x^k\pa_j\s\sum_i\left(\zD_{*k}^i\pa_i+\zD_{*k}^{*i}\xd
x^i\right)=-\zD^k_{*k}\pa_j+\zD_{*k}^{*j}\xd x^k.\ee The identity
(\ref{1.r1}) implies that $\zD^i_j=\zd^i_j\zD^k_k$ and
$-\zD^{*i}_j=\zd^i_k\zD^{*j}_k$. Since the indices $i,j,k$ are
arbitrary, we conclude that $\zD^i_j=\zl\zd^i_j$ for some
$\zl\in\R$ and $\zD^{*i}_j=0$, i.e., $\zD(\pa_j)=\zl\pa_j$.
Similarly, from the identity (\ref{1.r2}) we conclude that
$\zD(\xd x^j)=\zl'\xd x^j$. But now $\zl=\zl'$ follows from
(\ref{10}). Indeed, $(X+\zx)\s(X+\zx)=\xd i_X\zx$, so that $(\zl
X+\zl'\zx)\s(X+\zx)=\zl\pounds_X\zx-\zl'i_X\xd\zx=
(\zl-\zl')i_X\xd\zx+\zl\xd i_X\zx$ equals $\zl'\xd i_X\zx$ for all
vector fields $X$ and all 1-forms $\zx$, thus $\zl=\zl'$. \epf

\medskip\noindent
\textbf{Definition 3\ } A $(1,1)$-tensor on a Courant algebroid we
call \textit{orthogonal} if $N+N^*=0$. A (weak) Nijenhuis tensor
$N$ which is compatible with the symmetric pairing
$\la{\cdot},{\cdot}\ran$ of the Courant algebroid, in the sense
that (\ref{4}) and (\ref{5}) are satisfied for "$\s_N$" and
$\zr_N$, we call a \textit{(weak) Courant-Nijenhuis tensor}.

\medskip\noindent
Thus weak Courant-Nijenhuis tensors give rise to contractions of
Courant algebroids. Note however, that the structure of a Courant
algebroid is extremely rigid and that there are very few true
Courant-Nijenhuis tensors. First, observe that $N$ is a
Courant-Nijenhuis tensor if and only if $N-\frac{\zl}{2}I$ is
Courant-Nijenhuis (cf. \cite[Theorem 8]{CGM}), so we can always
reduce paired tensors to the case when $N+N^*=0$, i.e. to the case
of orthogonal $N$. Second, we have the following.
\begin{theo}\label{N}(\cite{CGMc}) If $N$ is an orthogonal Courant-Nijenhuis tensor,
then
$$X\s N^2Y=N^2(X\s Y),\quad \text{and}\quad N^2(Y\s Y)=(N^2Y)\s Y.$$
\end{theo}
\bepf
Using $N^*=-N$ and the invariance of the pairing, we get
\begin{equation}\label{p1}\la N(X\s_NY),Z\ran=-\la X\s_NY,NZ\ran=
-\zr(NX)\la Y,NZ\ran+\la Y,X\s_NNZ\ran
\end{equation}
and
\begin{equation}\label{p2}\la NX\s NY,Z\ran=\zr(NX)\la NY,Z\ran+\la Y,N(NX\s Z)\ran,
\end{equation}
so $N$ is Nijenhuis implies that the r.h. sides of (\ref{p1}) and
(\ref{p2}) are equal, i.e.
\begin{equation}\label{p3}
X\s_NNZ-N(NX\s Z)=0.
\end{equation}
But the l.h.s of (\ref{p3}) is
$$NX\s NZ-N(X\s_NZ)-N^2(X\s Z)+X\s N^2Z$$
and vanishing of the Nijenhuis torsion implies $N^2(X\s Z)=X\s
N^2Z$. The second identity one proves analogously, see the proof
of (\ref{10}). \epf
\begin{cor}\label{c1} Any Courant-Nijenhuis tensor $N$ on an
irreducible Courant algebroid satisfies:
\begin{description}
\item{(a)} \ $N+N^*=\zl I$,  \item{(b)} \ $N^2-\zl N+\zg I=0$,
\end{description}
for certain $\zl,\zg\in\R$, so that the algebra with involution
generated by $N$, thus the corresponding hierarchy, is trivial.
\end{cor}
\bepf According to Theorem \ref{t1.1}, $N+N^*=\zl I$. Then,
applying Theorem \ref{N} to $N:=N-\frac{\zl}{2}I$, we get
$(N-\frac{\zl}{2}I)^2=\zl' I$ which yields (b) with
$\zg=\frac{\zl^2}{4}-\zl'$. \epf

\bigskip\noindent
\textbf{Definition 4\ } An orthogonal Courant-Nijenhuis tensor $N$
on a Courant algebroid we call (for the terminology see \cite{BC})
\begin{description}\item{(i)} a \textit{complex Courant structure},
if $N^2=-I$; \item{(ii)} a \textit{product Courant structure}, if
$N^2=I$; \item{(iii)} a \textit{tangent Courant structure}, if
$N^2=0$.
\end{description}

\medskip\noindent
\textbf{Remark.} Note that complex Courant structures on the
canonical Courant algebroid $\cT M$ from Example 1 have been
introduced by N.~ Hitchin \cite{Hi} under the name of
\textit{generalized complex geometries}. They have been then
studied by M.~Gualtieri \cite{Gu} and have drawn an attention of
other authors (see e.g. \cite{Cr,LMTZ,Zu1,Zu2}). One can say, not
very precisely, that a \textit{generalized geometry} is a geometry
of contractions in which we replace a Nijenhuis tensor on the
tangent bundle (with the standard bracket of vector fields) with a
similar Nijenhuis tensor on the `extended tangent bundle' (with
the Courant bracket). When generalizing this scheme to an
arbitrary Courant algebroid, we can speak about a \textit{Courant
geometry}.

\begin{cor}\label{c2} Any orthogonal Courant-Nijenhuis tensor
on an irreducible Courant algebroid is proportional to either a
complex Courant structure, or to a product Courant structure, or
to a tangent Courant structure.
\end{cor}

\section{Courant geometries as supergeometries}
There is another approach to Courant algebroids, proposed by
D.~Roytenberg \cite{Ro,Ro2} (cf. also \cite{Vo}), in which the
Courant algebroid corresponds to a symplectic $\cN$-manifold
$(\wt{A},\zW)$ of degree 2 with the associated (graded) Poisson
bracket $\{\cdot,\cdot\}$, equipped additionally with a cubic
Hamiltonian $\zY$ which is \textit{homological}, i.e.
$\{\zY,\zY\}=0$. The symplectic $\cN$-manifold $\wt{A}$ of degree
2 is here the pullback of $\sT^*[2]A[1]$ (fibered canonically over
$(A\op A^*)[1]$) with respect to the embedding $A\hookrightarrow
A\op A^*$ given by $X\mapsto(X,\la X/2,\cdot\ran)$, i.e. it
completes the commutative diagram
\[
\begin{array}{ccc}
\wt{A} & \longrightarrow  & \sT^{*}[2]A[1]\\
\downarrow  &  & \downarrow \\
A[1] & \longrightarrow  & (A\oplus A^{*})[1]
\end{array}\]
Here we use the standard convention and, for a graded vector
bundle $E$ over a graded manifold $\cM$, write $E[n]$ for the
graded manifold obtained by shifting the fibre degrees by $n$. In
this picture, the corresponding Leibniz bracket is a
\textit{derived bracket} (cf. \cite{KS2,KS1}) for which $\zY$ is a
{\it generating Hamiltonian}:
\be\label{2.1} X\s Y=\{\{ X,\zY\},Y\}.
\ee
We should have probably written "$\s_\zY$" for the operation, but
let us fix $\zY$ and keep writing simply "$\s$". Note that the
above formula implies immediately
\be\label{2.2}\zr(X)(f)=\{\{
X,\zY\},f\}.
\ee
Here $X,Y$ are functions on $\wt{A}$ of degree 1 (i.e. sections of
$A$) and $f$ is of degree 0 (i.e. $f$ is a function on $M$).
Writing the graded algebra of super-functions on $\wt{A}$ as
$\A=\bigoplus_{k=0}^\infty\A^k$, we can identify the algebra of
functions on $M$ with $\A^0$ and the $\A^0$-module of sections of
$A$ with $\A^1$. The Poisson bracket reduced to $\A^1$ is just the
pseudo-Riemannian form $\la\cdot,\cdot\ran$. Moreover, the
Hamiltonian vector field $\pa_\zY=\{\zY,\cdot\}$ is a cohomology
operator in $\A$ defining the corresponding cohomology. The
pseudo-Riemannian form $\la\cdot,\cdot\ran$, thus the Poisson
bracket, identifies canonically $A$ with $A^*$ by $X\mapsto\{
X,\cdot\}$ (on sections) and any orthogonal $(1,1)$-tensor
$N\in\Sec(A\ot A^*)$ can be clearly identified with an element in
$\A^1\!\cdot\!\A^1\subset\A^2$, denoted, with some abuse of
notation, also by $N$. In this language, $N(X)=\{ N,X\}$ for
$X\in\A^1$. In an affine Darboux chart $(x^i,\zx^a,p_j)$ on
$\wt{A}$, corresponding to a chart $(x^i)$ on $M$ and a local
basis $\{ e_a\}$ of sections of $A$ such that $\la
e_a,e_b\ran=g_{ab}$, the symplectic form $\zW$ reads
$$\zW=\sum_i\xd p_i\xd x^i+\frac{1}{2}\sum_a\xd\zx^ag_{ab}\xd\zx^b$$
and $\zY\in\A^3$ is of the form
$$\zY=\sum_{a,i}\zx^a\zr^i_a(x)p_i-\frac{1}{6}\sum_{a,b,c}
\zvf_{a,b,c}(x)\zx^a\zx^b\zx^c,$$ where $\zr^i_a=\zr(e_a)(x^i)$
and $\zvf_{a,b,c}=\la e_a\s e_b,e_c\ran$. Any $(1,1)$-tensor
$N:A\ra A$, $N(e_a)=\sum_bN_a^b(x)e_b$ is orthogonal if and only
if $\sum_b\left(N_a^bg_{bc}+N_c^bg_{ab}\right)=0$ and then it is
represented by the element
$$N=\frac{1}{2}\sum_bN_a^bg_{bc}\zx^c\zx^a\in\A^1\cdot\A^1.$$

\medskip\noindent
Note, however, that what we have denoted $N^2$ before, and which
is $\{ N,\{ N,\cdot\}\}$ in the present notation, is not the
square of $N$ in the algebra $\A$. Since we will not use powers in
the algebra $\A$, we will keep the old notation.

\begin{prop} The derived bracket (\ref{2.1}) generated by any
cubic Hamiltonian $\zY$ always satisfies the compatibility
conditions (\ref{9}) and (\ref{10}). The Jacobi identity (\ref{3})
is equivalent to the homological condition $\{\zY,\zY\}=0$.
\end{prop}
\bepf  Since $\{ X,\{ Y,Y\}\}=0$ for $X,Y\in\A^1$, we have, due
to the graded Jacobi identity,
\beas 0&=&\{\zY,\{ X,\{
Y,Y\}\}\}=\{\{\zY, X\},\{ Y,Y\}\}\}- 2\{ X,\{\{\zY,
Y\},Y\}\}\\
&=&\zr(X)\la Y,Y\ran-2\la X,Y\s Y\ran,\eeas whence (\ref{9}). On
the other hand,
\beas\la X\s Y,Y\ran&=&\{\{\{\zY,X\},Y\},Y\}=\{\{\zY,X\},\{ Y,Y\}\}-
\{\{\{\zY,X\},Y\},Y\}\\
&=&\zr(X)\la Y,Y\ran-\la X\s Y,Y\ran,
\eeas
that proves (\ref{10}). That the Jacobi identity is equivalent to
the homological condition $\{\zY,\zY\}=0$ follows now from
\cite{Ro2}, Thorem 4.5.\epf

\medskip\noindent
\textbf{Example 2} (\cite{Ro2}) The symplectic $\cN$-manifold of
degree 2 associated with the canonical Courant algebroid from
Example 1 is
$$\wt{A}=\sT^*[2]\sT[1]M\simeq\sT^*[2]\sT^*[1]M$$ with the
canonical symplectic form $\zW$ of degree 2. In local affine
Darboux coordinates $(x^i,\zx^j,p_k,\zy_l)$, where $(x^i)$ are
local coordinates (of degree 0) on $M$, $(\zx^j,\zy_l)$ are
degree-1 coordinates associated with adapted linear functions on
the bundle $A=\cT M=\sT M\op\sT^*M$ corresponding to $\xd x^j$ and
$\pa_{x^l}$, and $(p_l)$ are degree-2 coordinates associated with
linear functions on another copy of $\sT^*M$ -- the core of the
double vector bundle $\sT^*[2]\sT[1]M\simeq\sT^*[2]\sT^*[1]M$.  In
these coordinates $\zW=\sum_i\left(\xd x^i\xd
p_i+\xd\zx^i\xd\zy_i\right)$ and the corresponding Poisson
superbracket reads
$$\{ F,G\}=i_{\xd G}i_{\xd F}\sum_i\left(\pa_{p_i}\pa_{x^i}+
\pa_{\zy_i}\pa_{\zx^i}\right).$$
The canonical cubic Hamiltonian
is in this case $\zY=\sum_i\zx^ip_i$ which is just the Hamiltonian
lift of the de Rham vector field $\xd=\sum_i\zx^i\pa_{x^i}$ on
$\sT[1]M$.

\begin{prop} For $N\in\A^1\!\cdot\!\A^1\subset\A^2$, the contracted
product $\s_N$ is the derived bracket associated with the cubic
Hamiltonian $\{\zY,N\}$. Moreover, $N$ represents a weak
Courant-Nijenhuis tensor if and only if $\{\{\zY,N\},N\}$ is a
$\pa_\zY$-cocycle, and $N$ represents a Courant-Nijenhuis tensor
if and only if $\{\{\zY,N\},N\}$ is the generating Hamiltonian for
"$\s_{N^2}$".
\end{prop}
\bepf We have
\beas\{\{ X,\{\zY,N\}\},Y\}&=&-\{\{ N,\zY\},X\},Y\}=\\
&=&-\{ N,\{\{\zY,X\},Y\}\}+\{\{\zY,\{ N,X\}\},Y\} +\{\{\zY,X\},\{
N,Y\}\}=\\
&=&-N(X\s Y)+NX\s Y+X\s NY=X\s_N Y.
\eeas
Thus, the product "$\s_N$" defines another Courant algebroid
structure on $(A,\la\cdot,\cdot\ran)$ if and only if
$\{\{\zY,N\},\{\zY,N\}\}=0$. But
$$\{\{\zY,N\},\{\zY,N\}\}=\{\{\{\zY,N\},\zY\},N\}\}\}-
\{\zY,\{\{\zY,N\},N\}\}=0-\{\zY,\{\{\zY,N\},N\}\}.$$ Since, as
easily seen, $X(\s_N)_NY=2\text{Tor}_N(X,Y)+X\s_{N^2}Y$, the
vanishing of the Nijenhuis torsion is equivalent that the
generator of "$(\s_N)_N$", i.e. $\{\{\zY,N\},N\}$ is the generator
of "$\s_{N^2}$". \epf
\begin{cor} A complex (resp., product, tangent) Courant structure
on a Courant algebroid $A$ associated with the qubic Hamiltonian
$\zY$ is exactly an element $N\in\A^1\cdot\A^1$ satisfying
$\{\{\zY,N\},N\}=-\zY$ (resp, $\{\{\zY,N\},N\}=\zY$,
$\{\{\zY,N\},N\}=0$).
\end{cor}
A detailed description of such quadratic Hamiltonians is in
general difficult, since in particular it contains all true
complex structures (cf. also \cite{Cr}). There are also relations
to presymplectic-Nijenhuis and Poisson-Nijenhuis structures, thus
bihamiltonian systems (cf. \cite[Theorem 10]{CGMc}).

\end{document}